\documentclass[12pt]{amsart}
\usepackage{amsthm, amssymb, amsfonts, mathrsfs, amsmath}

\usepackage{pb-diagram}
\usepackage{graphicx}
\usepackage{charter}
\usepackage[T1, T5]{fontenc}
\usepackage[a4paper,top=46mm,bottom=20mm,left=31mm,right=31mm]{geometry}

\usepackage{hyperref}
\newcommand\rurl[1]{%
  \href{http://#1}{\nolinkurl{#1}}%
}

\newtheorem*{theorem}{Main Theorem}

\numberwithin{equation}{section}
\numberwithin{definition}{section}
\numberwithin{lemma}{section}
\numberwithin{remark}{section}
\numberwithin{corollary}{section}
\numberwithin{proposition}{section}

\numberwithin{equation}{section}

\title[\uppercase{Non-trivial elements in the cohomology groups of the Steenrod algebra}]{\uppercase{A note on the non-trivial elements in the cohomology groups of the Steenrod algebra}}

\author[\DJ. V. Ph\'uc]{\DJ\d \abreve ng V\~o Ph\'uc}
\address{
Faculty of Education Studies\newline \indent
University of Khanh Hoa, 01 Nguyen Chanh, Khanh Hoa, Viet Nam}
\email{dangvophuc@ukh.edu.vn}



\begin{document}



\vspace*{-10mm}



\vspace{-2.8mm}


\vspace{-8.7mm}


\vspace{-7mm}

\vspace{1mm}

\begin{abstract}
Let $F_2$ be the prime field of two elements and let $GL_s:= GL(s, F_2)$ be the general linear group of rank $s.$ Denote by $\mathscr A$ the Steenrod algebra over $F_2.$ The (mod-2) Lambda algebra, $\Lambda,$ is one of the tools to describe those mysterious "Ext-groups". In addition, the $s$-th algebraic transfer of William Singer \cite{Singer} is also expected to be a useful tool in the study of them. This transfer is a homomorphism $Tr_s: F_2 \otimes_{GL_s}P_{\mathscr A}(H_{*}(B\mathbb V_s, F_2))\to {\rm Ext}_{\mathscr {A}}^{s,s+*}(F_2, F_2),$ where $\mathbb V_s$ denotes the elementary abelian $2$-group of rank $s$, and $H_*(B\mathbb V_s)$ is the homology group of the classifying space of $\mathbb V_s,$ while $P_{\mathscr A}(H_{*}(B\mathbb V_s, F_2))$ means the primitive part of $H_{*}(B\mathbb V_s, F_2)$ under the action of $\mathscr A.$ It has been shown that $Tr_s$ is highly non-trivial and, more precisely, that $Tr_s$ is an isomorphism for $s\leq 3.$ In addition, Singer proved that $Tr_4$ is an isomorphism in some internal degrees. He was also investigated the image of the fifth transfer by using invariant theory. 

In this note, we use another method to study the image of $Tr_5.$ More precisely, by direct computations using a representation of $Tr_5$ over the algebra $\Lambda,$ we show that $Tr_5$ detects the non-zero elements $h_0d_0\in {\rm Ext}_{\mathscr A}^{5, 5+14}(F_2, F_2),\ h_2e_0\in {\rm Ext}_{\mathscr A}^{5, 5+20}(F_2, F_2)$ and $h_1h_4c_0\in {\rm Ext}_{\mathscr A}^{5, 5+24}(F_2, F_2).$ The same argument can be used for homological degrees $s\geq 6$ under certain conditions.
\end{abstract}

\maketitle

\section{The Singer algebraic transfer and the (mod-2) Lambda algebra}

Let $\mathbb V_s$ be a rank $s$ elementary abelian 2-group and let $B\mathbb V_s$ denote the classifying space of $\mathbb V_s.$ We may equally well view $\mathbb V_s$ as an $s$-dimensional $F_2$-vector space. As well-known, $H^*(\mathbb V_s, F_2)\cong S(\mathbb V_s^{*}),$ the symmetric algebra over the dual space $\mathbb V_s^{*} = H^{1}(\mathbb V_s, F_2).$ Pick $x_1, \ldots, x_s$ to be a basis of $H^{1}(\mathbb V_s, F_2).$ Then, we have $H^{*}(B\mathbb V_s, F_2) = H^{*}(\mathbb V_s, F_2)\cong \mathcal P_s:= F_2[x_1, \ldots, x_s],$ the polynomial algebra on $s$ generators of degree one. It is considered as an unstable $\mathscr A$-module. The canonical $\mathscr A$-action on $\mathcal P_1$ is extended to an $\mathscr A$-action on $F_2[x_1, x_1^{-1}],$ the ring of finite Laurent series (see \cite{Adams2}, \cite{Wilkerson}). Then, $\overline{\mathcal P} = \langle \{x_1^{t}|\ t\geq -1\}\rangle$ is $\mathscr A$-submodule of $F_2[x_1, x_1^{-1}].$ One has a short-exact sequence:  
\begin{equation}\label{dkn1}
 0\to \mathcal P_1\xrightarrow{q} \overline{\mathcal P}\xrightarrow{\pi} \Sigma^{-1}F_2,
\end{equation}
where $q$ is the inclusion and $\pi$ is given by $\pi(x_1^{t}) = 0$ if $t\neq -1$ and $\pi(x_1^{-1}) = 1.$ Denote by $e_1$ the element in ${\rm Ext}_{\mathscr A}^1(\Sigma^{-1}F_2, \mathcal P_1),$ which is represented by the cocyle associated to \eqref{dkn1}. For each $i\geq 1,$ the short-exact sequence
\begin{equation}\label{dkn2}
 0\to \mathcal P_{i+1}\cong \mathcal P_i\otimes_{F_2}\mathcal P_1\xrightarrow{1\otimes_{F_2}q} \mathcal P_i\otimes_{F_2}\overline{\mathcal P}\xrightarrow{1\otimes_{F_2}\pi} \Sigma^{-1}\mathcal P_i,
\end{equation}
determines a class $(e_1\times \mathcal P_i)\in {\rm Ext}_{\mathscr A}^1(\Sigma^{-1}\mathcal P_i, \mathcal P_{i+1}).$ Then, using the cross product and Yoneda product, we have the element
$$ e_s = (e_1\times \mathcal P_{s-1})\circ (e_1\times \mathcal P_{s-2})\circ \cdots \circ (e_1\times \mathcal P_{1})\circ e_1 \in {\rm Ext}_{\mathscr A}^s(\Sigma^{-s}F_2, \mathcal P_s).$$
Let $\Delta(e_1\times \mathcal P_i): {\rm Tor}_{s-i}^{\mathscr A}(F_2, \Sigma^{-s}\mathcal P_i)\to {\rm Tor}_{s-i-1}^{\mathscr A}(F_2, \mathcal P_{i+1})$ be the connecting homomorphism associated to \eqref{dkn2}. Then, we have a composition of the connecting homomorphisms $$\overline{\varphi}_s = \Delta(e_1\times \mathcal P_{s-1})\circ \Delta(e_1\times \mathcal P_{s-2})\circ \cdots \circ \Delta(e_1\times \mathcal P_{1})\circ \Delta(e_1)$$from ${\rm Tor}_s^{\mathscr A}(F_2, \Sigma^{-s}F_2)$ to ${\rm Tor}_0^{\mathscr A}(F_2, \mathcal P_s) = F_2\otimes_{\mathscr A} \mathcal P_s$, determined by $\overline{\varphi}_s(z) = e_s\cap z$ for any $z\in {\rm Tor}_s^{\mathscr A}(F_2, \Sigma^{-s}F_2).$ Here $\cap$ denotes the \textit{cap product} in homology with $F_2$-coefficients. The image of $\overline{\varphi}_s$ is a submodule of the variant $(F_2 \otimes_{\mathscr A} \mathcal P_s)^{GL_s},$ where $GL_s$ denotes the general linear group of rank $s$ over $F_2.$ Hence, $\overline{\varphi}_s$ induces homomorphism
$$ \varphi_s: {\rm Tor}_s^{\mathscr A}(F_2, \Sigma^{-s}F_2)\to (F_2 \otimes_{\mathscr A} \mathcal P_s)^{GL_s}.$$
Since the supension $\Sigma^{-s}$ induces an isomorphism $${\rm Tor}_{s, *}^{\mathscr A}(F_2, \Sigma^{-s}F_2)\cong {\rm Tor}_{s, s+*}^{\mathscr A}(F_2, F_2),$$ we have the homomorphism
$$ \varphi_s: {\rm Tor}_{s, s+*}^{\mathscr A}(F_2, F_2)\to (F_2 \otimes_{\mathscr A} \mathcal P_s)_{*}^{GL_s}.$$
Let $P_{\mathscr A}(H_{*}(B\mathbb V_s, F_2))$ be the primitive part consisting of all elements in $H_{*}(B\mathbb V_s, F_2),$ which are annihilated by every positive degree operation in $\mathscr A.$ Then, the dual homomorphism
$$ Tr_s: F_2\otimes_{GL_s}P_{\mathscr A}(H_{*}(B\mathbb V_s, F_2))\to {\rm Ext}_{\mathscr A}^{s, s+*}(F_2, F_2)$$
of $\varphi_s$ is called \textit{the $s$-th Singer algebraic transfer}. It is expected to be a useful tool in describing the cohomology groups of $\mathscr A$ by means of invariant theory and the Peterson "hit problem" of determining a minimal set of $\mathscr A$-generatos for $\mathcal P_s$ (see \cite{Boardman}, \cite{Kameko}, \cite{Peterson}, \cite{Phuc-Sum1, Phuc-Sum2, Phuc3, Phuc4, Phuc5, Phuc6, Phuc8, Phuc9, Phuc10, Phuc11}, \cite{Walker-Wood}). It was then investigated by many authors (see Boardman \cite{Boardman}, Ch\ohorn n and H\`a \cite{Chon-Ha}, Crossley \cite{Crossley}, H\`a~\cite{Ha}, Minami \cite{Minami}, the present author \cite{Phuc3, Phuc4, Phuc6, Phuc7, Phuc8, Phuc10, Phuc11, Phuc12, Phuc13}, and others). By Boardman \cite{Boardman} and Singer \cite{Singer}, $Tr_s$ is an isomorphism for $s\leq 3.$ In \cite{Singer}, Singer also gave computations to show that $Tr_4$ is an isomorphism in a range of internal degrees. 

As it is known, the homological algebra $\{H_d(B\mathbb V_s, F_2)|\ d\geq 0\}$ is dual to $\mathcal P_s.$ Moreover, it is isomorphic to $\Gamma(a_1, \ldots, a_s),$ the divided power algebra generated by $a_1, \ldots, a_s,$ each of degree one, where $a_j = a_j^{(1)}$ is dual to $x_{j}.$ Here the duality is taken with respect to the basis of $\mathcal P_s$ consisting of all monomials in $x_1, \ldots, x_s.$ The algebra $\{H_d(B\mathbb V_s, F_2)|\ d\geq 0\}$ is a right $\mathscr A$-module. The right $\mathscr A$-action on this algebra is given by $(a^{(j)}_t)Sq^{s} = \binom{j-s}{s}a^{(j-s)}_t = Sq^{s}_*(a^{(j)}_t)$ and Cartan's formula. Note that $Sq^{s}_*$ is the dual Steenrod operation.

The (mod 2) Lambda algebra $\Lambda$ (see Bousfield et al. \cite{Bousfield}) is also one of  the tools to compute the cohomology groups of $\mathscr A.$ Recall that $\Lambda$ is defined as a differential, bigraded, associative algebra with unit over $F_2$ is generated by $\lambda_k\in\Lambda^{1,k},$ satisfying the Adem relations 
\begin{equation}\label{ct1}
\lambda_k\lambda_{2k+s+1} = \sum_{j\geq 0}\binom{s-j-1}{j}\lambda_{k+s-j}\lambda_{2k+1+j}\ \ (k\geq 0,\ s\geq 0)
\end{equation}
and differential 
\begin{equation}\label{ct2}
\partial (\lambda_{s-1}) = \sum_{j\geq 1}\binom{s-j-1}{j}\lambda_{s-j-1}\lambda_{j-1}\ \ (s\geq 1),
\end{equation}
where $\binom{s-j-1}{j}$ denotes the modulo $2$ value of the binomial coefficient, with the usual convention $\binom{s-j-1}{j} = 0$ for $j > s-j-1.$ In addition, $$H^{s, *}(\Lambda, F_2) = {\rm Ext}_{\mathscr{A}}^{s, s+*}(F_2, F_2).$$ Now, for non-negative integers $t_1, \ldots, t_s,$ a monomial  $\lambda_{t_1}\ldots\lambda_{t_s}\in \Lambda$ is called \textit{the monomial of the length $s$}. We shall write $\lambda_T,\ T = (t_1, \ldots, t_s)$ for $\prod_{1\leq k\leq s}\lambda_{t_k}$ and refer to $\ell(T) = s$ as the length of $T.$  Note that the algebra $\Lambda$ is not commutative and that the bigrading of a monomial indexed by $T$ may be written $(s, d),$ where the homological degree $s,$ as above, is the length of $T,$ and $d = \sum_{1\leq k\leq s}t_k.$ A monomial $\lambda_T$ is called \textit{admissible} if $t_k\leq 2t_{k+1}$ for all $1\leq k \leq s-1.$ By the relations \eqref{ct1}, the $F_2$-vector subspace $$\Lambda^{s, *} = \langle \{\lambda_T|\ T = (t_1, \ldots, t_s), \ t_k\geq 0, 1\leq k\leq s\}\rangle$$ of $\Lambda$ has an additive basis consisting of all admissible monomials of length $s.$ In \cite{Chon-Ha}, Ch\ohorn n and H\`a defined a homomorphism $\psi_s: H_*(B\mathbb V_s, F_2)\longrightarrow \Lambda^{s, *},$ which is considered as a representation in the algebra $\Lambda$ of the Singer transfer $Tr_s$ and determined by the following inductive formula:
$$ \psi_s(a^{T}) = \left\{ \begin{array}{ll}
\lambda_{t_1}&\text{if}\; \ell(T) = 1, \\
\sum_{j\geq t_s}\psi_{s-1}(\prod_{1\leq k\leq s-1}a_{k}^{(t_{k})}Sq^{j-t_s})\lambda_j &\text{if}\; \ell(T) > 1,
\end{array} \right.$$
for any $a^{T} := \prod_{1\leq k\leq s}a_{k}^{(t_{k})}\in H_*(B\mathbb V_s, F_2)$ and $T = (t_1, t_2, \ldots, t_{s}).$ Note that $\psi_s$ is not an $\mathscr A$-homomorphism. The authors showed in \cite{Chon-Ha} that if $\theta\in P_{\mathscr A}(H_*(B\mathbb V_s, F_2)),$ then $\psi_s(\theta)$ is a cycle in $\Lambda^{s, *}$ and $Tr_s([\theta]) = [\psi_s(\theta)]$. Therefrom, we have a homomorphism of algebras
$$ \psi = \{\psi_s:\ s\geqslant 0\}: \{H_*(B\mathbb V_s, F_2):\ s\geqslant 0\}\to \{\Lambda^{s, *}:\ s\geqslant 0\} =\Lambda,$$
which induces the Singer transfer.

By using invariant theory, Singer showed in \cite{Singer} that the non-zero element $Ph_1\in {\rm Ext}_{\mathscr A}^{5, 5+9}(F_2, F_2)$ was not detected by $Tr_5.$ The aim of the present note is also to study the image of $Tr_5.$ 

\section{Main result}
The main result of this note is the following.

\begin{theorem}\label{thmain} 
The non-zero elements
\begin{enumerate}
\item[(i)] $h_0d_0\in {\rm Ext}_{\mathscr A}^{5, 5+14}(F_2, F_2),$   

\item[(ii)] $h_2e_0 = h_0g_1\in {\rm Ext}_{\mathscr A}^{5, 5+20}(F_2, F_2),$

\item[(iii)] $h_1h_4c_0\in {\rm Ext}_{\mathscr A}^{5, 5+24}(F_2, F_2)$
\end{enumerate}
are in the image of the fifth algebraic transfer.
\end{theorem}


We prove the theorem by direct computations using the representation in the Lambda algebra of the fifth algebraic transfer. 

\begin{proof}[{\it Proof of main theorem}]

It is well-known that there exists an endomorphism $Sq^0$ of  the lambda algebra $\Lambda,$ determined by $Sq^0(\lambda_T = \prod_{1\leq k\leq s}\lambda_{t_k}) = \prod_{1\leq k\leq s}\lambda_{2t_k+1},$ where $\lambda_T$ is not necessarily admissible. It respects the Adem relations in \eqref{ct1} and commutes with the differential $\partial$ in \eqref{ct2} above. Then, $Sq^0$ induces the classical squaring operation in the Ext groups
$$Sq^0: H^{s, *}(\Lambda, F_2) = {\rm Ext}^{s, s+*}_{\mathscr {A}}(F_2, F_2) \to H^{s, s+2*}(\Lambda, F_2) = {\rm Ext}^{s, 2(s+*)}_{\mathscr {A}}(F_2, F_2).$$ According to Liulevicius \cite{Liulevicius}, $Sq^0$ is not the identity map. Moreover, it commutes with Kameko's squaring operation (see  \cite{Kameko}) $$ \widetilde {Sq^0}: F_2 \otimes_{GL_s}P_{\mathscr A}(H_{*}(B\mathbb V_s, F_2))\to F_2 \otimes_{GL_s}P_{\mathscr A}(H_{s+2*}(B\mathbb V_s, F_2))$$
through the $s$-th Singer transfer. In what follows, $$(Sq^{0})^{t}: {\rm Ext}^{*, *}_{\mathscr {A}}(F_2, F_2) \to {\rm Ext}^{*, *}_{\mathscr {A}}(F_2, F_2)$$ denotes the composite $Sq^{0}\ldots Sq^{0}$ ($t$ times of $Sq^{0}$)  if $t > 1,$ is $Sq^{0}$ if $t = 1,$ and  is the identity map if $t = 0.$ A family $\{a_t:\, t\geq 0\} \subset {\rm Ext}_{\mathscr A}^{s, s+*}(F_2, F_2)$ is called a $Sq^0$-family if $a_t = (Sq^0)^t(a_0)$ for $t\geq 0.$ According to Lin \cite{Lin}, the algebra $\{{\rm Ext}_{\mathscr A}^{s, s+*}(F_2, F_2):\ s\geqslant 0\}$ for $s\leqslant 4$ is generated by $h_t,\ c_t,\ d_t,\ e_t,\ f_t,\ g_{t+1},\ p_t,\ D_3(t),\ p'_t$ for $i\geqslant 0$ and subject to the relations $h_ih_{i+1} = 0,$ $h_ih^{2}_{i+2} = 0,$ and $h_i^3 = h^2_{i-1}h_{i+1}$ together with the relations  $h^2_th^2_{t+3} = 0, h_jc_t = 0$ for $j \in \{t-1,\ t,\ t+2,\, t+3\}.$   Furthermore, the set of the elements $d_t,\ e_t,\ f_t,\ g_{t+1},\ p_t,\ D_3(t)$ and $p'_t$ for $t\geqslant 0,$ is an $F_2$-basis for the indecomposable elements in ${\rm Ext}_{\mathscr A}^{4, 4+*}(F_2, F_2).$ By Lin \cite{Lin}, we have 
$$ \begin{array}{ll}
\medskip
h_t &= \{\lambda_{2^{t}-1} = (Sq^{0})^{t}(\lambda_0)\}\in {\rm Ext}_{\mathscr A}^{1, 2^{t}}(F_2, F_2) \ \mbox{for $t\geq 0$},\\
\medskip
c_t &= \{(Sq^{0})^{t}(\lambda_3^{2}\lambda_2)\}\in {\rm Ext}_{\mathscr A}^{3, 2^{t+3} + 2^{t+1}  +2^{t}}(F_2, F_2) \ \mbox{for $t\geq 0$},\\ 
\medskip
d_t &= \{(Sq^{0})^{t}(\lambda_3^2\lambda_2\lambda_6 + \lambda_3^2\lambda_4^2 + \lambda_3\lambda_5\lambda_4\lambda_2 + \lambda_7\lambda_1\lambda_5\lambda_1)\}\\
\medskip
&\quad\in {\rm Ext}_{\mathscr A}^{4, 2^{t+4} + 2^{t+1}}(F_2, F_2) \ \mbox{for $t\geq 0$},\\ 
e_t&= \{(Sq^{0})^{t}(\lambda_3^{3}\lambda_8 + (\lambda_3\lambda_5^{2} + \lambda_3^{2}\lambda_7)\lambda_4 + (\lambda_3^{2}\lambda_9 + \lambda_9\lambda_3^{2})\lambda_2)\}\\ 
\medskip
&\quad \in {\rm Ext}_{\mathscr A}^{4, 2^{t+4} + 2^{t+2}  +2^{t}}(F_2, F_2) \ \mbox{for $t\geq 0$},\\
g_{t+1}&= \{(Sq^{0})^{t}(\lambda_7^{2}\lambda_0\lambda_6 + (\lambda_3^{2}\lambda_9 + \lambda_9\lambda_3^{2})\lambda_5 +(\lambda_3^{2}\lambda_{11}+\lambda_3\lambda_9\lambda_5)\lambda_3)\}\\
&\quad\in {\rm Ext}_{\mathscr A}^{4, 2^{t+4} + 2^{t+3}}(F_2, F_2) \ \mbox{for $t\geq 0$}.
\end{array}$$
Combining this and the previous results by Tangora \cite{Tangora}, we get
$$ {\rm Ext}_{\mathscr A}^{5, 5+n}(F_2, F_2) = \left\{\begin{array}{ll}
\langle h_0d_0 \rangle &\mbox{if $n = 14$},\\
\langle h_2e_0 \rangle &\mbox{if $n = 20$},\\
\langle h_1h_4c_0 \rangle &\mbox{if $n = 24$},
\end{array}\right.$$
and $h_0d_0\neq 0,\ h_2e_0 = h_0g_1\neq 0,\ h_1h_4c_0 = h_3e_0\neq 0.$ Here, $h_t = [\lambda_{2^{t}-1}]$ is the Adams element (see \cite{Adams}) in ${\rm Ext}_{\mathcal{A}_2}^{1, 2^{t}}(F_2, F_2),$ for $0\leq t\leq 4.$ 

{\bf The case (i)}. Obviously, $\lambda_0\in \Lambda^{1, 0},$ and 
$$ \overline{d}_0 = \lambda_3^2\lambda_2\lambda_6 + \lambda_3^2\lambda_4^2 + \lambda_3\lambda_5\lambda_4\lambda_2 + \lambda_7\lambda_1\lambda_5\lambda_1 \in \Lambda^{4, 14}$$
 are the cycles in the algebra $\Lambda.$ Recall that $$h_0 = [\lambda_0]\in {\rm Ext}_{\mathscr{A}}^{1, 1}(F_2, F_2),\ \ d_0 = [\overline{d}_0]\in {\rm Ext}_{\mathscr{A}}^{4, 18}(F_2, F_2).$$ 
We consider the element $u_{14}$ in $H_{14}(B\mathbb V_5, F_2),$ which is the following sum:
$$ \begin{array}{ll}
&a_1^{(0)}a_2^{(6)}a_3^{(5)}a_4^{(2)}a_5^{(1)}+
a_1^{(0)}a_2^{(6)}a_3^{(4)}a_4^{(3)}a_5^{(1)}+
a_1^{(0)}a_2^{(6)}a_3^{(3)}a_4^{(4)}a_5^{(1)}+
\medskip
a_1^{(0)}a_2^{(6)}a_3^{(2)}a_4^{(5)}a_5^{(1)}\\
&\quad+
a_1^{(0)}a_2^{(6)}a_3^{(1)}a_4^{(6)}a_5^{(1)}+
a_1^{(0)}a_2^{(5)}a_3^{(6)}a_4^{(1)}a_5^{(2)}+
a_1^{(0)}a_2^{(5)}a_3^{(5)}a_4^{(2)}a_5^{(2)}+
\medskip
a_1^{(0)}a_2^{(5)}a_3^{(4)}a_4^{(3)}a_5^{(2)}\\
&\quad+
a_1^{(0)}a_2^{(5)}a_3^{(3)}a_4^{(4)}a_5^{(2)}+
a_1^{(0)}a_2^{(5)}a_3^{(2)}a_4^{(5)}a_5^{(2)}+
a_1^{(0)}a_2^{(5)}a_3^{(1)}a_4^{(6)}a_5^{(2)}+
\medskip
a_1^{(0)}a_2^{(5)}a_3^{(5)}a_4^{(1)}a_5^{(3)}\\
&\quad+
a_1^{(0)}a_2^{(3)}a_3^{(6)}a_4^{(2)}a_5^{(3)}+
a_1^{(0)}a_2^{(6)}a_3^{(2)}a_4^{(3)}a_5^{(3)}+
a_1^{(0)}a_2^{(6)}a_3^{(1)}a_4^{(4)}a_5^{(3)}+
\medskip
a_1^{(0)}a_2^{(5)}a_3^{(2)}a_4^{(4)}a_5^{(3)}\\
&\quad+
a_1^{(0)}a_2^{(3)}a_3^{(4)}a_4^{(4)}a_5^{(3)}+
a_1^{(0)}a_2^{(3)}a_3^{(2)}a_4^{(6)}a_5^{(3)}+
a_1^{(0)}a_2^{(3)}a_3^{(6)}a_4^{(1)}a_5^{(4)}+
\medskip
a_1^{(0)}a_2^{(3)}a_3^{(5)}a_4^{(2)}a_5^{(4)}\\
&\quad+
a_1^{(0)}a_2^{(3)}a_3^{(4)}a_4^{(3)}a_5^{(4)}+
a_1^{(0)}a_2^{(3)}a_3^{(3)}a_4^{(4)}a_5^{(4)}+
a_1^{(0)}a_2^{(3)}a_3^{(2)}a_4^{(5)}a_5^{(4)}+
\medskip
a_1^{(0)}a_2^{(3)}a_3^{(1)}a_4^{(6)}a_5^{(4)}\\
&\quad+
a_1^{(0)}a_2^{(5)}a_3^{(3)}a_4^{(1)}a_5^{(5)}+
a_1^{(0)}a_2^{(6)}a_3^{(1)}a_4^{(2)}a_5^{(5)}+
a_1^{(0)}a_2^{(5)}a_3^{(2)}a_4^{(2)}a_5^{(5)}+
\medskip
a_1^{(0)}a_2^{(3)}a_3^{(4)}a_4^{(2)}a_5^{(5)}\\
&\quad+
a_1^{(0)}a_2^{(5)}a_3^{(1)}a_4^{(3)}a_5^{(5)}+
a_1^{(0)}a_2^{(3)}a_3^{(3)}a_4^{(3)}a_5^{(5)}+
a_1^{(0)}a_2^{(3)}a_3^{(1)}a_4^{(5)}a_5^{(5)}+
\medskip
a_1^{(0)}a_2^{(6)}a_3^{(1)}a_4^{(1)}a_5^{(6)}\\
&\quad+
a_1^{(0)}a_2^{(5)}a_3^{(2)}a_4^{(1)}a_5^{(6)}+
a_1^{(0)}a_2^{(3)}a_3^{(4)}a_4^{(1)}a_5^{(6)}+
a_1^{(0)}a_2^{(3)}a_3^{(3)}a_4^{(2)}a_5^{(6)}+
a_1^{(0)}a_2^{(6)}a_3^{(6)}a_4^{(1)}a_5^{(1)}.
\end{array}.$$
Then, it is $\mathscr A^{+}$-annihilated. Indeed, by the unstable condition, we need to compute the actions of the Steenrod squares $Sq^{1}, Sq^{2}$ and $Sq^{4}.$ It is not difficult to check that $(u_{14})Sq^{i} = 0$ for $i = 1, 2.$ By direct calculation, $(u_{14})Sq^{4}$ is the following sum:
$$ \begin{array}{ll}
&a_1^{(0)}a_2^{(5)}a_3^{(3)}a_4^{(1)}a_5^{(1)}+
a_1^{(0)}a_2^{(3)}a_3^{(5)}a_4^{(1)}a_5^{(1)}+
a_1^{(0)}a_2^{(5)}a_3^{(3)}a_4^{(1)}a_5^{(1)}+
\medskip
a_1^{(0)}a_2^{(3)}a_3^{(5)}a_4^{(1)}a_5^{(1)}\\
&\quad+
a_1^{(0)}a_2^{(3)}a_3^{(3)}a_4^{(3)}a_5^{(1)}+
a_1^{(0)}a_2^{(3)}a_3^{(3)}a_4^{(3)}a_5^{(1)}+
a_1^{(0)}a_2^{(5)}a_3^{(1)}a_4^{(3)}a_5^{(1)}+
\medskip
a_1^{(0)}a_2^{(3)}a_3^{(1)}a_4^{(5)}a_5^{(1)}\\
&\quad+
a_1^{(0)}a_2^{(5)}a_3^{(1)}a_4^{(3)}a_5^{(1)}+
a_1^{(0)}a_2^{(3)}a_3^{(1)}a_4^{(5)}a_5^{(1)}+
a_1^{(0)}a_2^{(3)}a_3^{(5)}a_4^{(1)}a_5^{(1)}+
\medskip
a_1^{(0)}a_2^{(5)}a_3^{(3)}a_4^{(1)}a_5^{(1)}\\
&\quad+
a_1^{(0)}a_2^{(3)}a_3^{(3)}a_4^{(2)}a_5^{(2)}+
a_1^{(0)}a_2^{(5)}a_3^{(3)}a_4^{(1)}a_5^{(1)}+
a_1^{(0)}a_2^{(3)}a_3^{(5)}a_4^{(1)}a_5^{(1)}+
\medskip
a_1^{(0)}a_2^{(3)}a_3^{(2)}a_4^{(3)}a_5^{(2)}\\
&\quad+
a_1^{(0)}a_2^{(3)}a_3^{(3)}a_4^{(3)}a_5^{(1)}+
a_1^{(0)}a_2^{(3)}a_3^{(3)}a_4^{(2)}a_5^{(2)}+
a_1^{(0)}a_2^{(3)}a_3^{(3)}a_4^{(3)}a_5^{(1)}+
\medskip
a_1^{(0)}a_2^{(3)}a_3^{(2)}a_4^{(3)}a_5^{(2)}\\
&\quad+
a_1^{(0)}a_2^{(5)}a_3^{(1)}a_4^{(3)}a_5^{(1)}+
a_1^{(0)}a_2^{(3)}a_3^{(1)}a_4^{(5)}a_5^{(1)}+
a_1^{(0)}a_2^{(3)}a_3^{(1)}a_4^{(5)}a_5^{(1)}+
\medskip
a_1^{(0)}a_2^{(5)}a_3^{(1)}a_4^{(3)}a_5^{(1)}\\
&\quad+
a_1^{(0)}a_2^{(3)}a_3^{(3)}a_4^{(1)}a_5^{(3)}+
a_1^{(0)}a_2^{(3)}a_3^{(3)}a_4^{(1)}a_5^{(3)}+
a_1^{(0)}a_2^{(3)}a_3^{(1)}a_4^{(3)}a_5^{(3)}+
\medskip
a_1^{(0)}a_2^{(3)}a_3^{(1)}a_4^{(3)}a_5^{(3)}\\
&\quad+
a_1^{(0)}a_2^{(3)}a_3^{(2)}a_4^{(2)}a_5^{(3)}+
a_1^{(0)}a_2^{(3)}a_3^{(1)}a_4^{(3)}a_5^{(3)}+
a_1^{(0)}a_2^{(3)}a_3^{(2)}a_4^{(2)}a_5^{(3)}+
\medskip
a_1^{(0)}a_2^{(3)}a_3^{(1)}a_4^{(3)}a_5^{(3)}\\
&\quad+
a_1^{(0)}a_2^{(3)}a_3^{(3)}a_4^{(1)}a_5^{(3)}+
a_1^{(0)}a_2^{(3)}a_3^{(3)}a_4^{(2)}a_5^{(2)}+
a_1^{(0)}a_2^{(3)}a_3^{(3)}a_4^{(1)}a_5^{(3)}+
\medskip
a_1^{(0)}a_2^{(3)}a_3^{(2)}a_4^{(3)}a_5^{(2)}\\
&\quad+
a_1^{(0)}a_2^{(5)}a_3^{(1)}a_4^{(3)}a_5^{(1)}+
a_1^{(0)}a_2^{(3)}a_3^{(1)}a_4^{(5)}a_5^{(1)}+
a_1^{(0)}a_2^{(3)}a_3^{(2)}a_4^{(3)}a_5^{(2)}+
\medskip
a_1^{(0)}a_2^{(3)}a_3^{(1)}a_4^{(3)}a_5^{(3)}\\
&\quad+
a_1^{(0)}a_2^{(3)}a_3^{(1)}a_4^{(3)}a_5^{(3)}+
a_1^{(0)}a_2^{(3)}a_3^{(3)}a_4^{(1)}a_5^{(3)}+
a_1^{(0)}a_2^{(5)}a_3^{(1)}a_4^{(1)}a_5^{(3)}+
\medskip
a_1^{(0)}a_2^{(3)}a_3^{(1)}a_4^{(1)}a_5^{(5)}\\
&\quad+
a_1^{(0)}a_2^{(3)}a_3^{(2)}a_4^{(2)}a_5^{(3)}+
a_1^{(0)}a_2^{(5)}a_3^{(1)}a_4^{(1)}a_5^{(3)}+
a_1^{(0)}a_2^{(3)}a_3^{(1)}a_4^{(1)}a_5^{(5)}+
\medskip
a_1^{(0)}a_2^{(3)}a_3^{(2)}a_4^{(2)}a_5^{(3)}\\
&\quad+
a_1^{(0)}a_2^{(3)}a_3^{(3)}a_4^{(1)}a_5^{(3)}+
a_1^{(0)}a_2^{(3)}a_3^{(1)}a_4^{(3)}a_5^{(3)}+
a_1^{(0)}a_2^{(3)}a_3^{(1)}a_4^{(3)}a_5^{(3)}+
\medskip
a_1^{(0)}a_2^{(5)}a_3^{(1)}a_4^{(1)}a_5^{(3)}\\
&\quad+
a_1^{(0)}a_2^{(3)}a_3^{(1)}a_4^{(1)}a_5^{(5)}+
a_1^{(0)}a_2^{(3)}a_3^{(1)}a_4^{(1)}a_5^{(5)}+
a_1^{(0)}a_2^{(5)}a_3^{(1)}a_4^{(1)}a_5^{(3)}+
\medskip
a_1^{(0)}a_2^{(3)}a_3^{(3)}a_4^{(1)}a_5^{(3)}\\
&\quad+
a_1^{(0)}a_2^{(3)}a_3^{(3)}a_4^{(1)}a_5^{(3)}.
\end{array}$$
Thus, $(u_{14})Sq^{4} = 0.$ Then, using the differential \eqref{ct2} and the representation of $Tr_5$ over $\Lambda,$ we see that\\[1mm]
$\begin{array}{ll}
\psi_5(u_{14}) &= \lambda_0\lambda_3^2\lambda_2\lambda_6 + \lambda_0\lambda_3^2\lambda_4^2 + \lambda_0\lambda_3\lambda_5\lambda_4\lambda_2 + \lambda_0\lambda_7\lambda_1\lambda_5\lambda_1 + 
\partial(\lambda_0\lambda_3^2\lambda_9 + \lambda_0\lambda_3\lambda_9\lambda_3)\\
&= \lambda_0\overline{d}_0 + \partial(\lambda_0\lambda_3^2\lambda_9 + \lambda_0\lambda_3\lambda_9\lambda_3).
\end{array}$\\[1mm]
is a cycle in $\Lambda^{5, 14}$ and is a representative of $h_0d_0\in {\rm Ext}_{\mathscr A}^{5, 5+14}(F_2, F_2).$ Therefore, $h_0d_0\in {\rm Im}(Tr_5).$ 

\medskip

{\bf The case (ii)}. Obviously, $\lambda_3\in \Lambda^{1, 3}$ and $$ \overline{e}_0 = \lambda_3^{3}\lambda_8 + (\lambda_3\lambda_5^{2} + \lambda_3^{2}\lambda_7)\lambda_4 + \lambda_7\lambda_5\lambda_3\lambda_2 + \lambda_3^{2}\lambda_5\lambda_6\in \Lambda^{4, 17}$$ are the cycles in the algebra $\Lambda$ and they are representative of the non-zero elements $h_2\in {\rm Ext}_{\mathscr{A}}^{1, 4}(F_2, F_2)$ and $e_0\in {\rm Ext}_{\mathscr{A}}^{4, 21}(F_2, F_2)$ respectively. We consider the element $u_{20}\in H_{20}(B\mathbb V_5, F_2)$, which is the following sum:
\newpage
$$  \begin{array}{ll}
&
 a_1^{(3)}a_2^{(5)}a_3^{(5)}a_4^{(5)}a_5^{(2)}+
 a_1^{(3)}a_2^{(5)}a_3^{(5)}a_4^{(6)}a_5^{(1)}+
 a_1^{(3)}a_2^{(3)}a_3^{(5)}a_4^{(8)}a_5^{(1)}+
\medskip
 a_1^{(3)}a_2^{(5)}a_3^{(3)}a_4^{(8)}a_5^{(1)}\\
&\quad +
 a_1^{(3)}a_2^{(3)}a_3^{(6)}a_4^{(7)}a_5^{(1)}+
 a_1^{(3)}a_2^{(5)}a_3^{(7)}a_4^{(4)}a_5^{(1)}+
a_1^{(3)}a_2^{(7)}a_3^{(5)}a_4^{(4)}a_5^{(1)} + 
\medskip
 a_1^{(3)}a_2^{(3)}a_3^{(9)}a_4^{(4)}a_5^{(1)}\\
&\quad +
 a_1^{(3)}a_2^{(9)}a_3^{(3)}a_4^{(4)}a_5^{(1)}+
 a_1^{(3)}a_2^{(3)}a_3^{(9)}a_4^{(3)}a_5^{(2)}+
 a_1^{(3)}a_2^{(9)}a_3^{(3)}a_4^{(3)}a_5^{(2)}+
\medskip
 a_1^{(3)}a_2^{(5)}a_3^{(9)}a_4^{(2)}a_5^{(1)}\\
&\quad +
 a_1^{(3)}a_2^{(9)}a_3^{(5)}a_4^{(2)}a_5^{(1)} +
 a_1^{(3)}a_2^{(5)}a_3^{(10)}a_4^{(1)}a_5^{(1)}+
 a_1^{(3)}a_2^{(9)}a_3^{(6)}a_4^{(1)}a_5^{(1)}+
\medskip
 a_1^{(3)}a_2^{(3)}a_3^{(11)}a_4^{(2)}a_5^{(1)}\\
&\quad +
 a_1^{(3)}a_2^{(11)}a_3^{(3)}a_4^{(2)}a_5^{(1)} +
 a_1^{(3)}a_2^{(5)}a_3^{(5)}a_4^{(3)}a_5^{(4)}+
 a_1^{(3)}a_2^{(5)}a_3^{(3)}a_4^{(5)}a_5^{(4)}+
\medskip
 a_1^{(3)}a_2^{(3)}a_3^{(5)}a_4^{(5)}a_5^{(4)}\\
&\quad +
 a_1^{(3)}a_2^{(3)}a_3^{(12)}a_4^{(1)}a_5^{(1)}+
 a_1^{(3)}a_2^{(11)}a_3^{(4)}a_4^{(1)}a_5^{(1)}+
 a_1^{(3)}a_2^{(7)}a_3^{(8)}a_4^{(1)}a_5^{(1)}+
\medskip
 a_1^{(3)}a_2^{(7)}a_3^{(7)}a_4^{(1)}a_5^{(2)}\\
&\quad +  
a_1^{(3)}a_2^{(13)}a_3^{(2)}a_4^{(1)}a_5^{(1)} +
 a_1^{(3)}a_2^{(14)}a_3^{(1)}a_4^{(1)}a_5^{(1)}+
 a_1^{(3)}a_2^{(6)}a_3^{(5)}a_4^{(3)}a_5^{(3)}+
\medskip
 a_1^{(3)}a_2^{(5)}a_3^{(3)}a_4^{(6)}a_5^{(3)}\\
&\quad +
 a_1^{(3)}a_2^{(3)}a_3^{(6)}a_4^{(5)}a_5^{(3)} + 
a_1^{(3)}a_2^{(6)}a_3^{(3)}a_4^{(3)}a_5^{(5)}+
 a_1^{(3)}a_2^{(3)}a_3^{(3)}a_4^{(6)}a_5^{(5)}+
\medskip
 a_1^{(3)}a_2^{(3)}a_3^{(6)}a_4^{(3)}a_5^{(5)}\\
&\quad +
a_1^{(3)}a_2^{(5)}a_3^{(3)}a_4^{(3)}a_5^{(6)}+
 a_1^{(3)}a_2^{(3)}a_3^{(5)}a_4^{(3)}a_5^{(6)}+
 a_1^{(3)}a_2^{(3)}a_3^{(3)}a_4^{(5)}a_5^{(6)}+
\medskip
 a_1^{(3)}a_2^{(3)}a_3^{(3)}a_4^{(3)}a_5^{(8)}\\
&\quad + 
a_1^{(3)}a_2^{(3)}a_3^{(3)}a_4^{(4)}a_5^{(7)} +
 a_1^{(3)}a_2^{(3)}a_3^{(5)}a_4^{(2)}a_5^{(7)}+
 a_1^{(3)}a_2^{(3)}a_3^{(6)}a_4^{(1)}a_5^{(7)}+
\medskip
 a_1^{(3)}a_2^{(3)}a_3^{(3)}a_4^{(9)}a_5^{(2)}\\
&\quad +
 a_1^{(3)}a_2^{(3)}a_3^{(3)}a_4^{(10)}a_5^{(1)} + 
a_1^{(3)}a_2^{(5)}a_3^{(3)}a_4^{(7)}a_5^{(2)}+
 a_1^{(3)}a_2^{(5)}a_3^{(7)}a_4^{(3)}a_5^{(2)}+
 a_1^{(3)}a_2^{(7)}a_3^{(5)}a_4^{(3)}a_5^{(2)}.
\end{array}$$
Then it is $\mathscr A^{+}$-annihilated. Indeed, it is not difficult to verify that $(u_{20})Sq^{i} = 0$ for $i = 1, 2.$ A direct calculation shows that $(u_{20})Sq^{4}$ is the following sum:\\[1mm]
$\begin{array}{ll}
& 
a_1^{(3)}a_2^{(3)}a_3^{(5)}a_4^{(3)}a_5^{(2)}+
a_1^{(3)}a_2^{(5)}a_3^{(3)}a_4^{(3)}a_5^{(2)}+
a_1^{(3)}a_2^{(3)}a_3^{(3)}a_4^{(6)}a_5^{(1)}+
\medskip
a_1^{(3)}a_2^{(3)}a_3^{(5)}a_4^{(4)}a_5^{(1)}\\
&\quad+
a_1^{(3)}a_2^{(3)}a_3^{(3)}a_4^{(6)}a_5^{(1)}+
a_1^{(3)}a_2^{(5)}a_3^{(3)}a_4^{(4)}a_5^{(1)}+
a_1^{(3)}a_2^{(3)}a_3^{(3)}a_4^{(6)}a_5^{(1)}+
\medskip
a_1^{(3)}a_2^{(3)}a_3^{(7)}a_4^{(2)}a_5^{(1)}\\
&\quad+
a_1^{(3)}a_2^{(7)}a_3^{(3)}a_4^{(2)}a_5^{(1)}+
a_1^{(3)}a_2^{(3)}a_3^{(5)}a_4^{(4)}a_5^{(1)}+
a_1^{(3)}a_2^{(3)}a_3^{(7)}a_4^{(2)}a_5^{(1)}+
\medskip
a_1^{(3)}a_2^{(5)}a_3^{(3)}a_4^{(4)}a_5^{(1)}\\
&\quad+
a_1^{(3)}a_2^{(7)}a_3^{(3)}a_4^{(2)}a_5^{(1)}+
a_1^{(3)}a_2^{(3)}a_3^{(5)}a_4^{(3)}a_5^{(2)}+
a_1^{(3)}a_2^{(5)}a_3^{(3)}a_4^{(3)}a_5^{(2)}+
\medskip
a_1^{(3)}a_2^{(5)}a_3^{(5)}a_4^{(2)}a_5^{(1)}\\
&\quad+
a_1^{(3)}a_2^{(3)}a_3^{(7)}a_4^{(2)}a_5^{(1)}+
a_1^{(3)}a_2^{(5)}a_3^{(5)}a_4^{(2)}a_5^{(1)}+
a_1^{(3)}a_2^{(7)}a_3^{(3)}a_4^{(2)}a_5^{(1)}+
\medskip
a_1^{(3)}a_2^{(5)}a_3^{(6)}a_4^{(1)}a_5^{(1)}\\
&\quad+
a_1^{(3)}a_2^{(5)}a_3^{(6)}a_4^{(1)}a_5^{(1)}+
a_1^{(3)}a_2^{(3)}a_3^{(7)}a_4^{(2)}a_5^{(1)}+
a_1^{(3)}a_2^{(7)}a_3^{(3)}a_4^{(2)}a_5^{(1)}+
\medskip
a_1^{(3)}a_2^{(3)}a_3^{(3)}a_4^{(3)}a_5^{(4)}\\
&\quad+
a_1^{(3)}a_2^{(3)}a_3^{(5)}a_4^{(3)}a_5^{(2)}+
a_1^{(3)}a_2^{(5)}a_3^{(3)}a_4^{(3)}a_5^{(2)}+
a_1^{(3)}a_2^{(3)}a_3^{(3)}a_4^{(3)}a_5^{(4)}+
\medskip
a_1^{(3)}a_2^{(3)}a_3^{(3)}a_4^{(5)}a_5^{(2)}\\
&\quad+
a_1^{(3)}a_2^{(5)}a_3^{(3)}a_4^{(3)}a_5^{(2)}+
a_1^{(3)}a_2^{(3)}a_3^{(3)}a_4^{(3)}a_5^{(4)}+
a_1^{(3)}a_2^{(3)}a_3^{(3)}a_4^{(5)}a_5^{(2)}+
\medskip
a_1^{(3)}a_2^{(3)}a_3^{(5)}a_4^{(3)}a_5^{(2)}\\
&\quad+
a_1^{(3)}a_2^{(7)}a_3^{(4)}a_4^{(1)}a_5^{(1)}+
a_1^{(3)}a_2^{(7)}a_3^{(4)}a_4^{(1)}a_5^{(1)}+
a_1^{(3)}a_2^{(3)}a_3^{(3)}a_4^{(3)}a_5^{(4)}+
\medskip
a_1^{(3)}a_2^{(3)}a_3^{(3)}a_4^{(5)}a_5^{(2)}\\
&\quad+
a_1^{(3)}a_2^{(3)}a_3^{(3)}a_4^{(6)}a_5^{(1)}+
a_1^{(3)}a_2^{(3)}a_3^{(3)}a_4^{(5)}a_5^{(2)},
\end{array}$\\[1mm]
and  therefore by the unstable condition, $u_{20}\in P_{\mathscr A}(H_{20}(B\mathbb V_5, F_2)).$ Then, using the differential \eqref{ct2} the representation of the rank 5 transfer over the lambda algebra we see that\\[1mm]
$ \begin{array}{ll}
\psi_5(u_{20}) &=\lambda_3(\lambda_3^{3}\lambda_8 + (\lambda_3\lambda_5^{2} + \lambda_3^{2}\lambda_7)\lambda_4 + 
\lambda_7\lambda_5\lambda_3\lambda_2 + \lambda_3^{2}\lambda_5\lambda_6)\\
&\quad + \lambda_3^{2}\lambda_5\lambda_4\lambda_5 + \lambda_3^{2}\lambda_5\lambda_3\lambda_6 + 
\lambda_3^{2}\lambda_5\lambda_6\lambda_3\\
&= \lambda_3\overline{e}_0 + \partial(\lambda_3^{2}\lambda_5\lambda_{10} + \lambda_3^{2}\lambda_{12}\lambda_3 + \lambda_3\lambda_4\lambda_7^{2} + \lambda_3\lambda_0\lambda_{11}\lambda_7).
\end{array}$\\[1mm]
is a cycle in $\Lambda^{5, 20}$ and is a representative of the non-zero $h_2e_0\in {\rm Ext}_{\mathscr A}^{5, 5+20}(F_2, F_2).$ Thus,  $h_2e_0$ is in the image of $Tr_5.$

\medskip

{\bf The case (iii)}. We see that $\lambda_1\in \Lambda^{1, 1},\, \lambda_{15}\in \Lambda^{1, 15}$ and $\lambda_3^2\lambda_2\in \Lambda^{3, 8}$ are cycles in the algebra $\Lambda$ and that $h_1 = [\lambda_1]\in {\rm Ext}_{\mathscr{A}}^{1, 2}(F_2, F_2),\, h_4 = [\lambda_{15}]\in {\rm Ext}_{\mathscr{A}}^{1, 16}(F_2, F_2)$ and $c_0 = [\lambda_3^2\lambda_2]\in {\rm Ext}_{\mathscr{A}}^{3, 11}(F_2, F_2).$ We consider the following element in $H_{24}(B\mathbb V_5, F_2)$:
$$ 
u_{24} = 
 a_1^{(1)}a_2^{(15)}a_3^{(3)}a_4^{(3)}a_5^{(2)}+
 a_1^{(1)}a_2^{(15)}a_3^{(3)}a_4^{(4)}a_5^{(1)}+
 a_1^{(1)}a_2^{(15)}a_3^{(5)}a_4^{(2)}a_5^{(1)}+
 a_1^{(1)}a_2^{(15)}a_3^{(6)}a_4^{(1)}a_5^{(1)}.$$
Then, it is $\mathscr A^{+}$-annihilated. Indeed, by the unstable condition, we need only to consider the effects of the Steenrod squares $Sq^1$ and $Sq^2.$ A direct computations shows that
$$ \begin{array}{ll}
(u_{24})Sq^{1} &= a_1^{(1)}a_2^{(15)}a_3^{(3)}a_4^{(3)}a_5^{(1)} + a_1^{(1)}a_2^{(15)}a_3^{(3)}a_4^{(3)}a_5^{(1)}\\
\medskip
&\quad + a_1^{(1)}a_2^{(15)}a_3^{(5)}a_4^{(1)}a_5^{(1)} + a_1^{(1)}a_2^{(15)}a_3^{(5)}a_4^{(1)}a_5^{(1)}   = 0,\\
\medskip
(u_{24})Sq^{2} &=  a_1^{(1)}a_2^{(15)}a_3^{(3)}a_4^{(2)}a_5^{(1)} + a_1^{(1)}a_2^{(15)}a_3^{(3)}a_4^{(2)}a_5^{(1)}  = 0.
\end{array}$$
Now, using the representation in $\Lambda$ of the fifth algebraic transfer, we get
$$ \begin{array}{ll}
\medskip
\psi_5(a_1^{(1)}a_2^{(15)}a_3^{(3)}a_4^{(3)}a_5^{(2)}) &= \lambda_1\lambda_{15}(\lambda_3^{2}\lambda_2 + \lambda_3\lambda_4\lambda_1+  \lambda_4\lambda_3\lambda_1),\\
\medskip
\psi_5( a_1^{(1)}a_2^{(15)}a_3^{(3)}a_4^{(4)}a_5^{(1)}) &=\lambda_1\lambda_{15}(\lambda_3\lambda_4\lambda_1  + \lambda_4\lambda_3\lambda_1+\lambda_5\lambda_2\lambda_1),\\
\medskip
  \psi_5( a_1^{(1)}a_2^{(15)}a_3^{(5)}a_4^{(2)}a_5^{(1)}) &= \lambda_1\lambda_{15}(\lambda_5\lambda_2\lambda_1 + \lambda_6\lambda_1^{2}),\\
\medskip
 \psi_5(a_1^{(1)}a_2^{(15)}a_3^{(6)}a_4^{(1)}a_5^{(1)}) &=\lambda_1\lambda_{15}\lambda_1\lambda_6\lambda_1^{2},
\end{array}$$
 Then, since $u_{24}\in P_{\mathscr A}(H_{24}(B\mathbb V_5, F_2)),$  $\psi_5(u_{24}) = \lambda_1\lambda_{15}\lambda_3^{2}\lambda_2$ is a cycle in $\Lambda^{5, 24}$ and is a representative of the non-zero element $h_1h_4c_0\in {\rm Ext}_{\mathscr{A}}^{5, 29}(F_2, F_2).$ This implies that $h_1h_4c_0\in {\rm Im}(Tr_5).$ The proof of the main theorem is completed.

\end{proof}

\end{document}